\documentclass[11pt]{amsart}

\usepackage{amssymb,amsmath,enumerate,epsfig,rotating,graphics,changebar,eepic}
\usepackage{amsfonts, booktabs}
\usepackage{a4,latexsym,parskip}
\usepackage{hyperref}
\usepackage[all]{xy}

\hypersetup{pdfauthor=MS} \hypersetup{pdftitle=Fields}

\newcommand{\PP}{\mathbb{P}}

\newcommand{\C}[1][]{\mathbb{C}^{#1}}
\newcommand{\Q}[1][]{\mathbb{Q}_{#1}}
\newcommand{\N}[1][] {\mathbb{N}_{#1}}

\newcommand{\F}{\mathbb{F}}
\newcommand{\Z}{\mathbb{Z}}
\newcommand{\p}{\mathfrak{p}}

\newcommand{\OO}{\mathcal{O}}

\newcommand{\Aut}{\mathop{\rm Aut}\nolimits}
\newcommand{\NS}{\mathop{\rm NS}\nolimits}
\newcommand{\Gal}{\mathop{\rm Gal}\nolimits}
\newcommand{\disc}{\mathop{\rm disc}\nolimits}
\newcommand{\MW}{\mathop{\rm MW}\nolimits}
\newcommand{\MWL}{\mathop{\rm MWL}\nolimits}
\newcommand{\SL}{\mathop{\rm SL}\nolimits}
\newcommand{\corr}{\mathop{\rm corr}\nolimits}

\newcommand{\Qbar}{{\,\overline{\!\Q\!}\,}}

\def\NeS{N\'eron--Severi}
\def\MoW{Mordell--\mbox{\kern-.12em}Weil}
\def\SI{Shioda--Inose}
\def\0{^{\phantom0}}

\theoremstyle{break} 
\newtheorem{Theorem}{Theorem}

\newtheorem{Lemma}[Theorem]{Lemma}

\newtheorem{Algorithm}[Theorem]{Algorithm}

\theoremstyle{remark}
\newtheorem{Remark}[Theorem]{Remark}
\newtheorem{Example}[Theorem]{Example}

\theoremstyle{definition}

\begin{document}
\setlength{\unitlength}{1cm}

\title{Modular forms and K3 surfaces}


\author{Noam D.~Elkies}
\address{
Mathematics Department,
Harvard University,
1 Oxford Street,
Cambridge, MA 02138,
USA}
\email{elkies@math.harvard.edu}
\urladdr{http://www.math.harvard.edu/$\sim$elkies}

\author{Matthias Sch\"utt}
\address{Institut f\"{u}r Algebraische Geometrie, Leibniz Universit\"{a}t
  Hannover, Welfengarten 1, 30167 Hannover, Germany}
\email{schuett@math.uni-hannover.de}
\urladdr{http://www.iag.uni-hannover.de/$\sim$schuett/}






\date{March 22, 2013}

\begin{abstract}
For every known Hecke eigenform of weight 3
with rational eigenvalues we exhibit a K3 surface over $\Q$
associated to the form.  This answers a question asked independently by
Mazur and van Straten.  The proof builds on a classification of
CM forms by the second author.
\end{abstract}

\thanks{Partial funding from
the National Science Foundation (grant DMS-0501029)
and Deutsche Forschungsgemeinschaft
(research grants Schu 2266/2-1 and Schu 2266/2-2)
is gratefully acknowledged.}

\maketitle


\section{Introduction}
\label{s:intro}

The question of modularity for algebraic varieties over $\Q$
has been studied in great detail in recent years.
Historically, it began with work by A.~Weil on Fermat varieties \cite{W},
continued in the context of curves  
by Deuring  \cite{Deuring}
and  Eichler \cite{Eichler}.
Shimura then proved that every Hecke eigenform of weight~$2$
is associated to an abelian variety over~$\Q$ (conf.~\cite[\S 7]{Shimura}).
In the case of rational eigenvalues,
the corresponding variety is an elliptic curve.

Conversely, the Taniyama--Shimura--\mbox{\kern-.12em}Weil conjecture
states that every elliptic curve over $\Q$ is modular.
The celebrated proof of this conjecture by
Wiles et al.~\cite{Wi,TW,BCDT} not only implies Fermat's Last Theorem,
but also catalyzed many further developments in this area,
notably the proof of Serre's conjecture
by Khare and Wintenberger \cite{KW}.
This implies modularity for several classes of varieties,
for instance
rigid Calabi--Yau threefolds over $\Q$ (cf.~\cite{D}, \cite{GY}).
These results were preceded by work by Livn\'e \cite{L} on modularity for two-dimensional
motives with complex multiplication (CM) that we use here,
citing it as Theorem~\ref{Thm:mod}.

On the other hand, the problem of geometric realizations
is harder for Hecke eigenforms of weight greater than two.
Deligne \cite{Deligne} gives a geometric construction of
$\ell$-adic Galois representations for Hecke eigenforms.
However, the varieties involved vary greatly with the level.
In this sense, Deligne's construction is not as uniform as one might wish
(cf.~Remark~\ref{Rem:EMS}).

This paper solves the first case of higher weight where we can realize
{\bf all known Hecke eigenforms} with rational eigenvalues in a single class of varieties.
Conjecturally these comprise all Hecke eigenforms in question, as stated in our main theorem:

\begin{Theorem}\label{thm}
Assume the extended Riemann Hypothesis (ERH) for odd real Dirichlet
characters.  Then every Hecke eigenform of weight~$3$
with rational eigenvalues is associated to a K3 surface over $\Q$.
\end{Theorem}

This result answers a question asked independently by
Mazur and van Straten.  It builds on the classification of
CM forms with rational coefficients by the second author
which we recall in section \ref{s:CM}.
That section also explains the dependence of Theorem~\ref{thm} on the ERH.
Section \ref{s:K3} recalls the notion of singular K3 surfaces
and Livn\'{e}'s modularity result (Thm.~\ref{Thm:mod}).
We review the relevant known examples and
obstructions in sections \ref{s:Ex} and~\ref{s:obs}.
For every known Hecke eigenform of weight $3$ with rational eigenvalues
we then exhibit an explicit singular K3 surface over $\Q$,
thus proving Theorem \ref{thm}.
Our main technique to achieve this 
is constructing one-dimensional families of K3 surfaces
and searching for singular specializations over $\Q$.
This is explained in section \ref{s:tech}
and exhibited in detail for one particular family in section~\ref{s:fam}.
The paper concludes with the remaining surfaces
needed to prove Theorem~\ref{thm}.
A summary of the proof is given   in Section \ref{s:proof}.

\section{Singular K3 surfaces}
\label{s:K3}

A K3 surface is a smooth, projective, simply connected surface~$X$\/
with trivial canonical bundle $\omega_X=\OO_X$.
The most prominent examples are smooth quartics in $\PP^3$
and Kummer surfaces.  Later we will work with elliptic K3 surfaces.

Throughout this paper, modularity will refer to classical modular forms
(cf.~sect.~\ref{s:CM}).  This classical kind of modularity
is a very special property of a variety;
a general K3 surface over $\Q$ cannot be modular
for several reasons (cf.~the discussion before Theorem~\ref{Thm:mod}),
though the Langlands Program predicts a
correspondence with some automorphic forms.

K3 surfaces and their moduli have been studied in great detail.
We will come back to these questions in section \ref{s:tech}.
The only complex K3 surfaces that can be classically modular
are those that have no moduli at all.
In terms of the Picard number $\rho(X) = \text{rk}\NS(X)$,
the condition that $X$\/ have no moduli is that
\[
\rho(X) = 20,
\]
the maximum in characteristic zero.
K3 surfaces with Picard number $20$ are often referred to as
\emph{singular K3 surfaces}.
The terminology is reflected in the \SI\ structure
(cf.~sections~\ref{s:Ex}, \ref{s:obs})
which relates any singular K3 surface to a product
of two isogenous elliptic curves with complex multiplication~(CM),
thus with \emph{singular moduli}.

Our results will often be stated in terms of the \emph{discriminant} $d=d(X)$
of a singular K3 surface~$X$, i.e.~the discriminant of the
intersection form on the \NeS\ lattice,
which is the \NeS\ group endowed with the cup-product pairing:
\[
d=d(X)=\disc(\NS(X)).
\]

\begin{Example}\label{Ex:Fermat}
The Fermat quartic in $\PP^3$,
\[
S=\{(x_0 : x_1 : x_2 : x_3) \in \PP^3 \mid x_0^4+x_1^4+x_2^4+x_3^4=0\}
\]
is a singular K3 surface. There are $48$ obvious lines on~$S$,
all defined over the 8th cyclotomic field $\Q(e^{2\pi i/8})$.
Pjatecki\u\i -\v Sapiro and \v Safarevi\v c \cite{PSS} gave an argument
that the lines generate $\NS(S)$,
with $\disc(\NS(S)) = d = -64$. However, the argument depended
on a claim by Demjanenko whose proof turned out to be incorrect.
The proof was later independently completed by Cassels \cite{Cassels}
and Inose \cite{Inose-Kummer};
for an alternative approach see \cite{SSvL}.
\end{Example}

For any smooth projective surface $X$ over $\C$, we define the
\emph{transcendental lattice} $T(X)$
as the following sublattice of $H^2(X,\Z)$:
\begin{eqnarray}\label{eq:T}
T(X)=\NS(X)^\bot \subset H^2(X,\Z).
\end{eqnarray}
When $X$\/ is a K3 surface, $H^2(X,\Z)$ has rank $22$ and signature $(3,19)$.
Hence $T(X)$ has rank $22-\rho(X)$ and signature $(2,20-\rho(X))$.

If $X$ is a singular K3 surface, then $T(X)$ has rank~$2$, and is
even, positive-definite and equipped with an orientation
\cite[\S 4]{SI}.
Using the intersection form, we will identify the transcendental lattice
with a $2\times 2$ matrix
\[
T(X) \leftrightarrow \begin{pmatrix} 2a & b\\
b & 2c\end{pmatrix}
\]
with integral coefficients and discriminant $d=b^2-4ac<0$.
Applying an $\SL_2(\Z)$ change of basis,
we will always use a reduced representative with $-a < b\leq a\leq c$.
We will also use the shorthand $[2a,b,2c]$ for $T(X)$.

In consequence, if the singular K3 surface $X$ is defined
over some number field $L$,
then $T(X)$ gives rise to a two-dimensional Galois
representation $\varrho$ over $L$.
Over some extension of the ground field,
Shioda and Inose associated
a Hecke character to this Galois representation~\cite{SI}.
Livn\'e then proved modularity over $\Q$ as an application of
a more general result concerning motives:

\begin{Theorem}[Livn\'e {\cite{L}}]
\label{Thm:mod}
Every singular K3 surface over $\Q$ is modular:
there is a Hecke eigenform with system of
eigenvalues $\{c_p\}$, such that
\[
\text{trace $\varrho($Frob}_p) = c_p \;\;\; \text{for almost all } p.
\]
\end {Theorem}
The corresponding Hecke eigenform $f$ has weight $3$ and CM by
the imaginary quadratic field $K=\Q(\sqrt{d})$,
where $d=\disc(\NS(X))<0$ is the discriminant of~$X$.
The Hecke eigenvalues of $f$ are all rational integers.

For instance, the Hecke eigenform corresponding to the Fermat quartic
has CM by $\Q(\sqrt{-1})$.
For the model $S$ in Example~\ref{Ex:Fermat},
the associated newform is the $\eta$-product
$\eta(4\tau)^6$ of weight 3 and  level $16$
(see \cite[Table 2]{S-CM} for eigenvalues).

\section{CM newforms}
\label{s:CM}

Modular forms for congruence subgroups come with
weight, level and nebentypus character that appear in the
transformation law. For fixed invariants, the
modular forms constitute a module over the Hecke algebra.
On cusp forms, the Hecke operators
can be diagonalized simultaneously. The set of eigenvalues
as in Theorem~\ref{Thm:mod} determines a unique primitive
normalized eigenform, a so-called \emph{newform}.

By the transformation law, a modular form $f$ always has
a Fourier expansion
\[
f=f(\tau)=\sum_{n\geq 0} a_n q^n,\;\;\;q=e^{2\pi i \tau}.
\]
Here a (normalized) newform is characterized by the property that
the Fourier coefficients are multiplicative and 
$a_p$ equals the eigenvalue of the $p$-th Hecke operator.
In particular the Mellin transform (or $L$-series of $f$)
\[
L(f,s) = \sum_{n\geq 0} a_n n^{-s}
\]
has an Euler product. 
Throughout this paper, we will assume that any newform is normalized 
and use the terms Fourier coefficients
and Hecke eigenvalues interchangeably.

We are interested in newforms of weight $3$ with rational eigenvalues.
Generally, odd weight is special in the sense that the nebentypus
character is necessarily nontrivial. In fact, by a result of
Ribet \cite{R}, a newform $f$ of odd weight with real Hecke eigenvalues
has CM by its
nebentypus character
associated to some imaginary quadratic field $K$. In particular, $f$ comes from a Hecke character of $K$.

CM newforms with rational coefficients have been classified by the
second author in \cite{S-CM}. The analysis of the associated Hecke
characters revealed the following structure:

\begin{Theorem}[Sch\"utt {\cite{S-CM}}]
\label{Thm:CM}
For fixed weight $k+1$, there is a bijective correspondence
\[
\left\{\begin{matrix} \text{CM newforms of weight}\; k+1\\
\text{with rational eigenvalues}\\ \text{up to twisting}
       \end{matrix}\right\}
\stackrel{1:1}{\longleftrightarrow}\left\{\begin{matrix}
\text{Imaginary-quadratic fields $K$}\\\text{with class group }
Cl(K)\subseteq (\Z/k\Z)^g\\ \text{for some}\;\, g\in\N
                                          \end{matrix}\right\}
\]
\end{Theorem}

Unless $K=\Q(\sqrt{-1})$ or $\Q(\sqrt{-3})$, twisting refers to
modifying the Fourier coefficients by a quadratic Dirichlet
character $\chi$ (since otherwise the rationality of Fourier
coefficients is not preserved):
\begin{eqnarray}\label{eq:twist}
f\otimes \chi = \sum_{n\geq 1} a_n \chi(n) q^n.
\end{eqnarray}
Any twist of a newform is again a newform
(possibly after modifying some zero coefficients at primes 
dividing the conductor of the twisting character), but the level
and the nebentypus character will differ in general.
If $K=\Q(\sqrt{-1})$ or $\Q(\sqrt{-3})$, then we can also twist the
associated Hecke character by a biquadratic resp.~cubic character.
Twisting can often be achieved geometrically,
for instance for elliptic fibrations in Weierstrass form
(see the discussion around (\ref{eq:Kummer})
and in Section \ref{s:proof}).

For fixed weight, we will refer to the Hecke eigenforms
from Theorem~\ref{Thm:CM} as
newforms of class number $h(K)$.
In this terminology, it should be understood that the newforms
have rational Hecke eigenvalues.

Conjecturally, the number of imaginary quadratic fields
on the right hand side in Theorem~\ref{Thm:CM} is always
finite. The case of $k=1$ was proven by Heilbronn
and made explicit by Heilbronn and Linfoot \cite{Heil}.
The exponents $k=2, 3$ are due to Weinberger \cite{Wb}
(the latter also proven by Boyd and Kisilevsky \cite{BK}).
Recently, Heath-Brown proved finiteness for
$k=5$, $2^n$, or $3\cdot 2^n$ for $n\in\N$  \cite{H-B}.

Most of these results are not effective.
However, for $k=2$ (and $k=1$),
our knowledge goes much further.
There are 65 imaginary quadratic fields whose
class groups are at most two-torsion.
Their class numbers go up to $16$ for $K=\Q(\sqrt{-1365})$.
Many of these fields (or rather their discriminants) were already known to Euler
in 1778
through his search for idoneal or convenient numbers (\cite{euler}; see also \cite{Frei}).
We list them by discriminant $d_K$ and class number $h(K)$:

$$
\begin{array}{c|l}
\hline
h(K) & d_K\\[0.2ex]
\hline
\hline
1 & -3, -4, -7, -8, -11, -19, -43, -67, -163\\
\hline
2 & -15, -20, -24, -35, -40, -51, -52, -88, -91, -115, \\
& -123, -148, -187, -232, -235, -267, -403, -427\\
\hline
4 & -84, -120, -132, -168, -195, -228, -280, -312,\\
&   -340, -372, -408, -435, -483, -520, -532, -555,\\
& -595, -627, -708, -715, -760, -795, -1012, -1435\\
\hline
8 & -420, -660, -840, -1092, -1155, -1320, -1380,\\
& -1428, -1540, -1848, -1995, -3003, -3315\\
\hline
16 & -5460\\
\hline
\end{array}
$$

In \cite{Wb}, Weinberger proved that there is
at most one further imaginary quadratic field
with class group exponent~$2$.
Assuming the absence of Siegel--Landau zeros
for the \hbox{$L$-series} of odd real Dirichlet characters,
the known list is complete.
This condition would be a consequence of the extended
Riemann hypothesis (ERH) for odd real Dirichlet characters.
Thus the (stronger) assumption in Theorem~\ref{thm}.

This paper solves the geometric realization problem in weight~3.
By Theorem~\ref{Thm:CM}, our primary task is to find a
singular K3 surface over $\Q$ for each of the 65 newforms (up to twisting)
corresponding to imaginary quadratic fields with class group
exponent~$2$.
As we will work with elliptic surfaces,
we need not worry about twisting
(cf.~the discussion around (\ref{eq:twist}) and (\ref{eq:Kummer})
and in Section \ref{s:proof}).

\begin{Remark}[Elliptic modular surfaces]\label{Rem:EMS}
There is a canonical way to realize newforms of weight 3,
via the elliptic modular surface (or universal elliptic curve)
for the corresponding congruence subgroup $\Gamma_1(N)$
(see \cite{Deligne}, \cite{ShEMS}).
However, this correspondence is not uniform in several respects.
First, the complex surfaces involved vary greatly with $N$.
Even if one considers a newform and a twist as in (\ref{eq:twist}),
the corresponding surfaces will in general not be \hbox{$\Qbar$-isomorphic}.
Moreover, the surfaces will in general have more than
one associated newform (i.e.~$p_g = h^{2,0} > 1$).
\end{Remark}

There are two classes of surfaces which might come to
mind first when thinking of modular forms:
abelian surfaces and K3 surfaces ---
both generalizations of elliptic curves to dimension two.
They will be discussed in the next section.

\section{Singular abelian surfaces and Kummer surfaces}
\label{s:Ex}

Shioda and Inose derived a canonical way to produce a
singular K3 surface of given isomorphism class
(i.e.~of given transcendental lattice by the Torelli theorem) \cite{SI}.
The main object in their construction is the Kummer surface
Km$(A)$ of the abelian surface $A$ with the
given transcendental lattice.
Hence we shall briefly discuss singular abelian surfaces over $\Q$.

A complex abelian surface $A$ has $H^2(A, \Z)$ of rank~$6$.
The surface is called singular if $\rho(A)=4$,
that is, if its transcendental lattice $T(A)$ has rank~$2$.
By a result of Shioda and Mitani \cite{SM},
every singular abelian surface is isomorphic to a
product of two isogenous elliptic curves $E, E'$ with CM.
These can be given explicitly in terms of the transcendental lattice $T(A)$.

By the classical theory,
an elliptic curve with CM has a model over $\Q$ if and only if 
its endomorphism ring has class number one \cite{Shimura}.
Up to $\bar\Q$-isomorphism, there are exactly 13 such curves.
If both $E$ and $E'$ are defined over $\Q$,
we deduce that the product $A=E\times E'$ realizes a
newform of weight 3 with CM by the corresponding
imaginary quadratic field $K$ of class number~$1$.

For singular abelian surfaces over $\Q$,
there is one further possible construction:
If the imaginary quadratic field $K$ has class number~$2$,
then consider an elliptic curve $E$ with CM by the
ring of integers $\OO_K$ in $K$.
Then $E$ is defined over a quadratic extension of $\Q$.
Let $A$ be the Weil restriction of $E$ to $\Q$.
By \cite{SM}, $A$ is a singular abelian surface
with transcendental lattice corresponding to the
non-principal class in $Cl(K)$.
Since $E$ is associated to some Hecke character,
the abelian surface $A$ over $\Q$ geometrically realizes a newform of weight 3 with CM by $K$.

It follows that we can realize all newforms of weight 3 and
class number $1$ or~$2$ in singular abelian surfaces over $\Q$.
Obstructions from
the cohomological structure of singular abelian surfaces 
(beyond the analogues of Theorems \ref{Thm:genus} and \ref{Thm:NS})
prevent us from
pursuing this approach any further:

\begin{Lemma}\label{Lem:ab}
Let $A$ be a singular abelian surface of discriminant $d$.
Assume that $A$ is defined over $\Q$.
Then $\Q(\sqrt{d})$ has class number 1 or 2.
\end{Lemma}

\emph{Proof:}
Let $K=\Q(\sqrt{d})$.
It suffices to prove the lemma for fundamental discriminants $d$,
since otherwise analogous arguments apply to the even bigger Galois action
of the ring class field $H(d)$ over $\Q$.
By \cite[Thm.~5.4]{S-fields} (the abelian surface version of Theorem \ref{Thm:genus} to come),
the class group of $H(d)$ is only 2-torsion:
\[
\mbox{Cl}(d) = \Gal(H(d)/K) \cong (\Z/2\Z)^m.
\]
Equivalently the Hilbert class field $H(d)$ has a totally real subfield $F$ of degree $2$
which contains (and is generated by) the j-invariants of all elliptic curves with
CM by the ring of integers $\mathcal{O}_K$.
Here $F$ is Galois over $\Q$
with Galois group
\[
\Gal(F/\Q) \cong \Gal(H(d)/K) \cong (\Z/2\Z)^m.
\]
We shall now study the Galois action of $\sigma\in\Gal(F/\Q)$ on $\NS(A)$.
By \cite{SM}, $A$ is isomorphic over $\C$ to a product of elliptic curves with CM by $\mathcal O_K$:
\begin{eqnarray}
\label{eq:A}
A \cong_{\C} E \times E'.
\end{eqnarray}
Applying Galois elements $\sigma\in\Gal(F/\Q)$,
we derive equivalent decompositions
\[
A \cong_{\C} E^\sigma \times E'^\sigma.
\]
That is, $A$ contains any elliptic curve $E$ with CM by $\mathcal O_K$.
Conversely, an elliptic curve $E\subset A$ induces a decomposition as in \eqref{eq:A}
since the linear system $|E|$ induces an elliptic fibration $\pi$ with base another elliptic curve $E'$,
and the fibration is constant, that is all fibers are isomorphic to $E$ over $\C$.
Since $E^2=0$ by adjunction, we deduce the numerical inequivalence
\[
E \not\equiv E^\sigma \;\;  \text{ in } \NS(A) \;\;\; \forall \sigma\in\Gal(F/\Q),
\]
for otherwise $E.E^\sigma=0$ would imply that $E^\sigma$ is also a fiber of $\pi$
and therefore $\C$-isomorphic to $E$, but $j(E)\neq j(E^\sigma)=j(E)^\sigma$ by construction.

{\bf Conclusion:} Any non-trivial element $\sigma\in\Gal(F/\Q)$ acts non-trivially on $\NS(A)$.

Decomposing $\NS(A)\otimes\Q$ into the Galois-invariant class of a hyperplane section
and its orthogonal complement of rank three,
we infer readily that $m\leq 2$, i.e.~the class number is $h(d)=1,2,4$.
In order to rule out the third alternative and prove the lemma,
it suffices to show that also $\Gal(K/\Q)$ acts non-trivially on $\NS(A)$.
To this end, we relate the determinant of the Galois representation on $\NS(A)$
to the quadratic Dirichlet character associated to $K$:

{\bf Claim:} $\det \NS(A) = \chi_K$.

To prove the claim, we employ a similar reasoning as in \cite{Sh-corr} 
using that all cohomology of $A$ is determined by $H^1$ 
not only in singular cohomology,
but also in the $\ell$-adic \'etale case as Galois representations:
\[
H^2(A,{\Q}_\ell) = \wedge^2 H^1(A,{\Q}_\ell).
\]
Since $H^1(A,{\Q}_\ell)$ has determinant $\chi_\ell^2$ with the $\ell$-adic cyclotomic character
$\chi_\ell$,
we infer that $H^2(A,{\Q}_\ell)$ has determinant $\chi_\ell^6$.
Via the cycle class map, $\NS(A)$ embeds into $H^2(A,{\Q}_\ell(1))$
which also contains an orthogonal 2-dimensional transcendental piece $T$ coming from $T(A)$ as in \eqref{eq:T}.
Since $H^2(A,{\Q}_\ell(1))$ has determinant $1$,  
we read off that 
\[
\det \NS(A) = (\det T)^{-1}.
\]
The singular abelian surface $A$ over $\Q$ is modular;
in detail, by \cite{L}, the  newform attached to $T(A)$ has nebentypus character $\chi_K$.
This immediately translates to $\det T=\chi_K$ which proves the claim.

Summing things up, assuming that $h(d)=4$,
we have derived too big a Galois action to be accommodated by a rank 3 sublattice of $\NS(A)$.
This leaves only $h(d)=1,2$ as possibilities and thus proves Lemma \ref{Lem:ab}.
\qed

By Lemma \ref{Lem:ab}, we can only have singular abelian surfaces  over~$\Q$
for imaginary quadratic fields of class number $1$ and~$2$. In order to realize
the newforms of weight 3 and greater class number geometrically,
we therefore turn to singular K3 surfaces.
This has two advantages:
On the one hand, the above constructions for class number $1$
and~$2$ carry over to the corresponding Kummer surfaces
(cf.~Lemma \ref{Lem:Kummer}).
On the other, the cohomological obstructions are much less restrictive 
for singular K3 surfaces.
In detail, this will be studied in the next section.

We now turn to Kummer surfaces. Here we consider the quotient of
an abelian surface $A$ by its involution~$-1$. The minimal resolution of
the resulting 16 singularities of type $A_1$ is a K3 surface.
One multiplies the intersection form on $T(A)$ by~$2$
to obtain $T(\text{Km}(X))=T(A)[2]$.
For example, the Fermat surface (Example~\ref{Ex:Fermat})
is $\Qbar$-isomorphic to the Kummer surface
of  the product of elliptic curves with
CM by $\Z[\sqrt{-1}]$ and $\Z[2\sqrt{-1}]$.

\begin{Lemma}\label{Lem:Kummer}
Any newform of weight 3 and class number~$1$ or~$2$
is associated to a Kummer surface over $\Q$.
\end{Lemma}

\emph{Proof:}
The Kummer quotient is defined over the ground field of the abelian surface
(cf.~(\ref{eq:Kum})).
For each imaginary quadratic field $K$ of class number $1$ or~$2$,
we have found a singular abelian surface over $\Q$.
Therefore we obtain a Kummer surface over $\Q$ for some
newform with CM by $K$. We now address the issue of twisting.

We will use elliptic fibrations. Let $E, E'$ be elliptic curves.
Denote the $j$-invariants by $j, j'$.
Inose exhibited an explicit elliptic fibration on Km$(E\times E')$ in \cite{Inose}.
By \cite[proof of Prop.~4.1]{S-fields}, this elliptic fibration
can be defined over $L=\Q(j+j', j\cdot j')$:
\begin{eqnarray}\label{eq:Kum}
\text{Km}(E\times E'):\;\;\; y^2 = x^3 + B(t)\,x+C(t),\;\;\; B, C\in L[t].
\end{eqnarray}

In the present situation,
the abelian surface $A$ that we start with
is isomorphic to the product
of two isogenous CM elliptic curves of class number $1$ or~$2$.
Then $L=\Q$ because $j,j'$ are equal for class number~$1$
and quadratic conjugate for class number~$2$.
Let $f$ denote the associated newform.
Then the $\Q(\sqrt{d})$-isomorphic fibration
\begin{eqnarray}\label{eq:Kummer}
\text{Km}(A):\;\; d\,y^2 = x^3 + B(t) x + C(t), \;\;\; d\in\Q^*
\end{eqnarray}
realizes the twist of $f$ by the quadratic
Dirichlet character $\left(\frac d{\cdot}\right)$ as in (\ref{eq:twist}).
Hence all quadratic twists have geometric equivalents
(and this works generally for Weierstrass fibrations).
This leaves cubic and biquadratic twists in case of
$K=\Q(\sqrt{-3}), \Q(\sqrt{-1})$.
Here we use the fact that for the elliptic curves
with $j$-invariant, $j=0$ resp.~$1728$, the Kummer surface of the self product admits
an automorphism of order~$6$ resp.~$4$. In fact, in the first case,
Inose's elliptic fibration (\ref{eq:Kum}) is isotrivial ($B$ is the zero polynomial).
Hence every fiber admits an automorphism of order~$6$,
so we can apply cubic twists fiberwise.
In the second case we could also argue directly with the
Fermat quartic (Example~\ref{Ex:Fermat}), since it is Kummer. \qed

By Lemma \ref{Lem:Kummer}, we are left to find
singular K3 surfaces over $\Q$ for all newforms of weight 3
with class number 4, 8 or 16. Up to twisting, there are 38 such forms known.

\section{Obstructions for singular K3 surfaces}
\label{s:obs}

Before continuing our search for singular K3 surfaces over $\Q$,
we discuss obstructions to the field of definition.
These are much milder than for singular abelian surfaces,
notably because $H^1(X)$ is trivial for a K3 surface~$X$
(cf.~Lemma~\ref{Lem:ab}),
so cohomological obstructions are only imposed by $H^2(X)$.
Moreover, not every singular K3 surface is Kummer
although the relation with Kummer surfaces is very close.

Shioda and Inose showed in \cite{SI}
how to produce a singular K3 surface
with given transcendental lattice.
Every singular K3 surface $X$ admits a Nikulin involution
such that the quotient is Kummer.
The resulting picture is often referred to as
\emph{\SI\ structure}.
We sketch it in the following figure.
Here $A$ and $X$ are chosen with the same transcendental lattice
$T(X)=T(A)$, and $T($Km$(A))=T(A)[2]$.
 
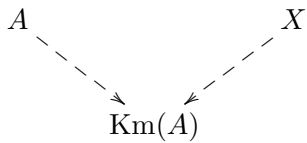
\[
  \xymatrix{A \ar@{-->}[dr] && X\ar@{-->}[dl]\\
 & \mbox{Km}(A)&}
 \]

The \SI\ construction is exhibited over some
finite extension of the field of definition of~$A$;
the elliptic fibration (\ref{eq:Kum}) on the Kummer surface
is a base change from~$X$
so that  Km$(A)$ is in fact sandwiched by $X$ as in \cite{Sandwich}.
The major question is when we can descend
$X$ to $\Q$. We now discuss the known obstructions.

The first obstruction comes from lattice theory.
It involves the genus of a lattice and was first studied
by Shimada \cite{Shimada} in the case of fundamental discriminant.
The second author then proved the general case in \cite{S-fields}.

\begin{Theorem}[Sch\"utt, Shimada]\label{Thm:genus}
Let $X$ be a singular K3 surface over some number field.
Let $d$ denote the discriminant of $X$ and $K=\Q(\sqrt{d})$. Then
\[
\{T(X^\sigma) \mid \sigma\in\Aut(\C/K)\} = \text{genus of } T(X).
\]
\end{Theorem}

If $X$ is defined over $\Q$, the genus of $T(X)$ consists of a single class.
This implies that $Cl(K)$ is at most two-torsion,
which of course also follows from Theorems~\ref{Thm:mod} and \ref{Thm:CM}.


The second obstruction is related to class field theory.
It essentially says that even if a singular K3 surface descends
to some number field $L$,
it still carries the structure of the ring class field $H(d)$
through the Galois action on the \NeS\ group.
This property was first noted for singular K3 surfaces over $\Q$
by the first author in \cite{Elkies}.
The second author pursued an alternative approach in \cite{S-NS}.

\begin{Theorem}[Elkies, Sch\"utt]\label{Thm:NS}
Let $X$ be a singular K3 surface of discriminant $d$.
Let $H(d)$ be the ring class field for $d$.
Let $L$ be a number field such that
$\NS(X)$ is generated by divisors over $L$.
Then $H(d)\subseteq L(\sqrt{d})$.
\end{Theorem}

\begin{Remark}
Conversely it was recently proved by Hulek and one of us
that any singular K3 surface $X$ of discriminant $d$
admts a model with $\NS(X)$ is generated by divisors over 
the ring class field $H(d)$ \cite{HS}.
\end{Remark}

In consequence, as $h(d)$ increases,
singular K3 surfaces become more and more complicated.
In particular the number of singular K3 surfaces over $\Q$
up to $\C$-isomorphism is finite.
This result is originally due to \v Safarevi\v c \cite{Shafa}.
Hence, in our attempt to find singular K3 surfaces for
all newforms of weight three with rational Fourier coefficients,
we are searching within a finite set.

We now resume our search for singular K3 surfaces over $\Q$
for the newforms of weight 3 and class numbers 4, 8 and 16.
Here we mention another class of singular K3
surfaces which has been investigated before:
\emph{extremal elliptic K3 surfaces}.
These are singular K3 surfaces admitting an elliptic fibration
with finite group of sections.
Up to the torsion sections,
such surfaces are determined by the configuration of singular fibers.
Extremal elliptic K3 fibrations were classified
by Shimada and Zhang in \cite{SZ}. They are finite in number.

Many explicit defining equations have been obtained
by the second author in \cite{S-Rocky} and
by Beukers and Montanus in \cite{BM}.
In addition to previous newforms of class number~$1$~and~$2$,
they realize ten discriminants of class number $4$ and $8$.
The next table lists the discriminants and one possible
configuration of singular fibers and \MoW\ group for each newform.
In each case the surface is semistable, that is, all reducible fibers
are of type $I_n$ for some $n\geq 1$;
so we simplify the notation by listing only the indices~$n$.
For completeness, we also give the transcendental
lattice $T_X$ in the shorthand notation
$[2a,b,2c]$ for the intersection form.

$$
\begin{array}{|c|c|c|c|}
\hline
\text{discriminant} & \text{configuration} & MW & T_X\\
\hline
\hline
-84 & [1,2,2,2,3,14] & \Z/2\Z & [2,0,42]\\
  \hline
 -120 & [1,2,2,2,5,12] & \Z/2\Z & [6,0,20]\\
   \hline
 -4\cdot 132 & [1,1,3,4,4,11] & (0) & [24, 12, 28]\\
   \hline
  -168 & [1,1,1,3,4,14] & (0) & [4, 0, 42]\\
  \hline
  -195 & [1,1,1,3,5,13] & (0) & [6, 3, 34]\\
    \hline
  -280 & [1,1,1,4,7,10] & (0) & [2, 0, 140]\\
    \hline
  -312 & [1,1,2,3,4,13] & (0) & [6, 0, 52]\\
  \hline
    \hline
  -4\cdot 420 & [1,3,4,4,5,7] & (0) & [24, 12, 76]\\
    \hline
  -660 & [1,2,2,3,5,11] & (0) & [2, 0, 330]\\
    \hline
  -840 & [1,1,4,5,6,7] & (0) & [12, 0, 70]\\
    \hline
\end{array}
$$

By Theorem \ref{Thm:NS}, the Galois group of the
ring class field $H(d)$ acts nontrivially
on the \NeS\ group $\NS(X)$,
i.e.~on the reducible fibers.
A nontrivial action can occur in two ways: 
either a reducible fiber is fixed under the Galois action, 
but it has type
$I_n$ for $n\geq 3$, $I_n^*$ for $n\geq 0$ or $IV, IV^*$, allowing for an involution of the components
(or for $I_0^*$ even an $S_3$ action)
that preserves both identity component and incidence relations;
or there are several reducible fibers of the same type lying over points on
the base curve $\PP^1$ that are Galois conjugates.
Later in this paper we again need a term for points on~$\PP^1$
lying under singular fibers (not necessarily reducible) of an elliptic surface;
we call them ``cusps''$\!$, consistent with the special case of
the universal elliptic curve over the modular curve $X(N)$ or $X_1(N)$.

In terms of the associated newforms,
the discriminants in the above table exhaust
the extremal elliptic K3 surfaces over $\Q$.
In the next section we will explain our main techniques
to exhibit singular K3 surfaces
for the remaining 28 imaginary quadratic fields.

\section{The main techniques}
\label{s:tech}

We search for singular K3 surfaces over $\Q$
with particular discriminants.
The main idea is to take advantage of the
moduli theory of complex K3 surfaces.
Any one-dimensional family of
K3 surfaces of Picard number $\rho\geq 19$ has
infinitely many specializations with $\rho=20$
(see for instance \cite{oguiso}).
We will only search for the specializations over $\Q$
(which are finite in number by \v Safarevi\v c's result \cite{Shafa}).
Because of Theorem~\ref{Thm:NS}, these surfaces
and therefore the families involved
will be very special.
Hence one of the key steps will be
to construct suitable families.
This will be achieved in the next section.
Here we explain how we find the specializations
of a given family.

Given a one-dimensional family of K3 surfaces $X_\lambda$ over $\Q$
satisfying $\rho(X_\lambda) \geq 19$ for all~$\lambda$,
there is an easily checked necessary condition that must be satisfied
by any member~$X$\/ with $\rho(X)=20$.
The condition is based on the Lefschetz fixed point formula
at a good prime $p$, formulated in terms of $\ell$-adic
\'etale cohomology $H_{\acute{e}t}^i(\bar X, {\Q}_\ell)$
for some prime $\ell\neq p$.
Here we work with the base change $\bar X$\/
of the reduction of $X \bmod p$
to an algebraic closure of $\F_p$. For simplicity, we will just write
$H^i(X)$ for $H_{\acute{e}t}^i(\bar X, {\Q}_\ell)$ in the following.

The cohomology groups $H^i(X)$ are equipped with an induced action of
the geometric Frobenius morphism Frob$_p$.
The set $X(\F_p)$ of $\F_p$-rational points
on $X$ is exactly the fixed set of Frob$_p$.
For a K3 surface $X/\Q$ with good reduction at~$p$,
the Lefschetz fixed point formula simplifies to
\begin{eqnarray}\label{eq:Lefschetz}
\# X(\F_p) = 1 + \text{trace Frob}_p^*(H^2(X)) + p^2.
\end{eqnarray}
Because $\rho(X)\geq 19$, we can predict 19 of the 22
eigenvalues of Frob$_p^*$ on $H^2(X)$.
Since the absolute Galois group operates through
a permutation on the algebraic cycles,
all these eigenvalues have the form $\zeta\cdot p$
for some roots of unity $\zeta$.
By the Weil conjectures,
one further eigenvalue has the form $\pm p$,
and the remaining two eigenvalues are algebraic
integers $\alpha_p,\beta_p$ of absolute value $p$.
In particular, the unordered pair $(\alpha_p, \beta_p)$
is determined by (\ref{eq:Lefschetz})
and the sign of the other eigenvalue $\pm p$.
In general, $\beta_p=\bar\alpha_p$;
the only exception is the supersingular case ($\rho(\bar X)=22$)
where $\alpha_p+\beta_p=0$.

If the specialization $X$ at some $\lambda_0\in\Q$
is a singular K3 surface,
then it is modular by Theorem~\ref{Thm:mod}.
Hence, for the right choice of sign in the third eigenvalue,
\begin{eqnarray}\label{eq:match}
\alpha_p+\beta_p=a_p
\end{eqnarray}
where $a_p$ is the Fourier coefficient of the
corresponding newform $f$ of weight 3.
Since $f$ has CM, both $\alpha_p$ and $\beta_p$ lie in the
imaginary quadratic extension $K$ associated to $f$.
By Theorem~\ref{Thm:CM}, $K$ has class group exponent $2$.
Moreover, $K$ remains fixed when $p$ varies.
This gives a criterion to either rule out $\lambda_0$ or
collect evidence for $\rho(X)=20$.

As it stands, our condition for $\rho=20$ is necessary but not sufficient.
To search for the CM-specializations,
we will use the condition in a different,
almost opposite approach:
we search for good parameters mod $p$ and
try to lift $\lambda_0$ to $\Q$.
We use the following algorithm:

\begin{Algorithm}\label{Alg}
Let $X_\lambda$ be a family of K3 surfaces over $\Q$
with Picard number $\rho\geq 19$.
Then the following algorithm returns candidate
parameters $\lambda_0\in\Q$ such that the
specialization $X$ at $\lambda_0$ might have $\rho(X)=20$:

\begin{enumerate}[(i)]
\item Fix one of the 65 weight 3 newforms $f$ up to twisting
(that is, one of the imaginary quadratic fields of class group exponent two).

\item \label{item}
For each of several primes~$p=p_i$ ($i=1,\ldots,n$),
  use  (\ref{eq:Lefschetz}) to compute $\alpha_p,\beta_p$
  for every $\lambda\in\F_p$, and find all $\lambda$ such that
  one choice of sign for the other eigenvalue leads to
  $\alpha_p, \beta_p$ matching $f$\/ as in (\ref{eq:match}).
  Even though $f$ varies, we need only compute $\#X_\lambda(\F_p)$
  once for each pair $(p,\lambda)$.

\item \label{item-3}
For a collection of parameters
$\lambda_1\bmod p_1, \hdots, \lambda_n \bmod p_n$
matching the newform $f$,
compute a lift $\lambda_0\in\Q$ of small height using the Chinese
Remainder Theorem and the Euclidean algorithm.

\end{enumerate}
\end{Algorithm}

For each newform $f$, the algorithm returns a
number of guesses for specializations in $\Q$, if any.
Often, there will be one $\lambda_0$ among them
which looks particularly likely
(small height, small primes involved etc.).
At this point, we can continue to collect numerical evidence
by running the above test for further primes $p$.
In the end, however, we want to prove that $\rho(X_{\lambda_0})=20$.
Therefore we have to find explicitly an additional divisor
on the specialization $X$ at $\lambda_0$.
This is where we turn to elliptic surfaces.
Until this point the procedure works for any one-dimensional family
(and generalizes to other settings such as Calabi-Yau threefolds
with third Betti number 4).
For instance, we computed candidates for the singular specializations
for the Dwork pencil
\[
x_0^4+x_1^4+x_2^4+x_3^4=\lambda x_0x_1x_2x_3
\]
of deformations of the Fermat quartic from~Example~\ref{Ex:Fermat}
(cf.~\cite{ES} for a detailed account).  But the only way we know
to systematically search for extra divisors on a K3 family $X_\lambda$
uses models of $X_\lambda$ as elliptic surfaces (with section).

The advantage of working with elliptic surfaces is the following.
By the formula of Shioda and Tate,
their \NeS\ groups are always generated by horizontal
and vertical divisors, i.e.~by fiber components and sections.
Hence for the Picard number to increase in a family,
either the singular fibers degenerate further
(which happens for only finitely many $\lambda_0 \in \C$),
or there is an additional section~$P$.
Then the discriminant $d$ of the specialization
can be computed purely in terms of the intersection behavior
of~$P$\/ with the singular fibers and the other sections.
This is made explicit through the theory of \MoW\ lattices, see \ref{ss:ell}.

In our setting, $d$ is predicted up to a square factor
by the newform $f$\/ and its CM-field $K$\/ by
Theorem~\ref{Thm:mod} and the explanation following it.
This provides us with additional information about
the conjectural section $P$, information that is often crucial
to the feasibility of a direct computation.
The next section illustrates this by the
detailed analysis of a particular family of K3 surfaces.

\section{Specializations of a one-dimensional family}
\label{s:fam}

We want to search in one-dimensional families of
K3 surfaces with Picard number $\rho\geq 19$
for singular specializations over $\Q$.
In fact, it is a nontrivial task to find such families
with interesting specializations in the first place.
This difficulty is due to the nontrivial Galois action
that the \NeS\ group must admit
by Theorem~\ref{Thm:NS}.
In this section, we discuss one particular family in detail.
In fact, we start out with a two-dimensional family.

\subsection{A two-dimensional family of elliptic K3 surfaces}

We start with the following
two-dimensional family of elliptic K3 surfaces
in extended Weierstrass form
with parameters $\lambda\in\PP^1, \mu\neq 0$:
\begin{eqnarray*}
X_{\lambda, \mu}:\;\;y^2 & = & x^3 +
(t-\lambda) A \,x^2 +  t^2 (t-1) (t-\lambda)^2 B\, x +
t^4 (t-1)^2 (t-\lambda)^3 C,\\
A & = & \frac 1{24} \Bigl( \frac 19 (2  \mu+9)^3 t^3-
(22  \mu-9) (2  \mu-27) t^2-27 (14  \mu-9) t-81 \Bigr),\\
B & = &  \mu\, \Bigl(\frac 19 (2  \mu+9)^3 t^2
-2 (10  \mu-9) (2  \mu-9) t-27 (2  \mu-3) \Bigr),\\
C & = &  \frac 23  \mu^2  ((2  \mu+9)^3 t-81 (2  \mu-3)^2).
\end{eqnarray*}
This elliptic surface has discriminant
$$
\Delta = 36  \mu^4 t^5 (t-1)^3 (t-\lambda)^6 h(t)
$$
where $h(t)=(2  \mu+9)^4 t^3-9 (32  \mu+27) (2  \mu+9)^2 t^2+
 81 (308  \mu^2+243-864  \mu) t+729 (4  \mu-9).$
 
It has the following singular fibers:
$$
\begin{array}{|c||c|c|c|c|c|}
\toprule
\text{cusp} & 0 & 1 & \infty & \lambda & \text{roots of } h\\
\midrule
\midrule
\text{fiber} & I_5 & I_3 & I_7 & I_0^* & I_1, I_1, I_1\\
\bottomrule
\end{array}
$$
In general, $X_{\lambda, \mu}$ has \NeS\ lattice
\[
\NS(X_{\lambda, \mu}) = U + A_2 + A_4 + A_6 + D_4.
\]
Here $U$ denotes the hyperbolic plane
(generated by the 0-section $O$\/ and a general fiber~$F$\/),
and $A_i, D_i$ are the root lattices
corresponding to the reducible singular fibers.
In particular, we deduce that
\[
\rho(X_{\lambda, \mu})\geq 18.
\]
We briefly explain how we found the above family.
As in \cite{E-Shimura},
we work with an extended Weierstrass form.
Here we can translate $x$
so that all singular fibers have their singularities at $x=y=0$.
In the present situation, this gives rise to the family
\begin{eqnarray*}
y^2 & = & x^3 +
(t-\lambda)\, a_2(t) \,x^2 +  t^2\, (t-1)\, (t-\lambda)^2\, a_1(t)\, x +
t^4\, (t-1)^2\, (t-\lambda)^3\, a_0(t).
\end{eqnarray*}
In general, the $a_i(t)$ are polynomials of degree $\deg(a_i)\leq i+1$.
Hence we have ten parameters to choose including $\lambda$,
relative to one normalization by rescaling $x, y$.
The above extended Weierstrass form guarantees
that the fiber types are at least $I_4$, $I_2$, $I_0^*$, $I_2$
at $t=0,1,\lambda,\infty$ respectively.
We can easily choose the coefficients of $a_2(t)$
to promote the fiber at $\infty$ to type $I_6$.
Then we solve a system of three nonlinear equations
in the five coefficients of $a_0(t), a_1(t)$
to derive the family $X_{\lambda, \mu}$.
This can be achieved by appropriate linear combinations
of the equations and a suitable choice of normalization.
Details are available from the authors upon request.

A few words on how we came to choose the family $X_{\lambda,\mu}$ might be in order.
To allow for a substantial Galois action, 
there ought to be  several reducible fibers.
For a singular K3 surface of big class number  over $\Q$ to appear as a  specialization of the family,
the fundamental discriminant eventually will have several distinct prime factors,
so some of these already have to occur in the discriminants of the root lattices associated to the reducible fibers.
Presently, these discriminants are, up to sign, $3, 5, 7, 4$
which shows that the family $X_{\lambda,\mu}$ is close to being optimal.
On the other hand, the family is still accessible to an explicit parametrization;
this is the main obstacle against working with families with even bigger discriminant
(which is why later we retreat to extended mod $p$ methods, see \ref{ss:5460}, \ref{ss:2}).

The family $X_{\lambda, \mu}$ can easily be specialized
to a family $X_\lambda$ with $\rho(X_\lambda)\geq 19$
by degenerating the singular fibers, i.e.~merging fibers.
There are two independent ways to do so.
On the one hand, we can match $\lambda$ with one of the cusps.
This results in four families.
Each of them has several interesting specializations.
The case where we merge $I_0^*$ with a $I_1$-fiber
will be taken up in the next section (Example~\ref{Ex:1}).
On the other hand, we can merge one of the $I_1$ fibers
with another singular fiber by a suitable choice of  $\mu$.
We now discuss the most beneficial case.

\subsection{A one-dimensional family of elliptic K3 surfaces}

Within the family of elliptic K3 surfaces $X_{\lambda,\mu}$,
we can merge two fibers of type $I_1$
by setting $\mu=\frac{243}{10}$.
We obtain a one-dimensional family of elliptic K3 surfaces $X_\lambda$
with the following singular fibers:
$$
\begin{array}{|c||c|c|c|c|c|c|}
\toprule
\text{cusp} & 0 & 1 & \frac{35}{32} & -\frac 5{1024} & \infty & \lambda\\
\midrule
\midrule
\text{fiber} & I_5 & I_3 & I_2 & I_1 & I_7 & I_0^*\\
\bottomrule
\end{array}
$$
The general member $X_\lambda$ has \NeS\ lattice
\[
\NS(X_\lambda) = U + A_1 + A_2 + A_4 + A_6 + D_4,
\]
so $\rho(X_\lambda)\geq 19$.
The Galois action on the \NeS\ group is encoded
in the fields where the singular fibers with at least three components split
(depending on $\lambda$):
$$
\begin{array}{|c||c|c|c|c|}
\toprule
\text{fiber} & I_5 & I_3 & I_7 & I_0^*\\
\midrule
\midrule
\text{splitting field} & \Q(\sqrt{6\lambda}) &
\Q(\sqrt{10\,(1-\lambda)}) & \Q(\sqrt{15}) & \Q(f(\lambda))\\
\bottomrule
\end{array}
$$
Here $\Q(f(\lambda))$ is the splitting field of the cubic polynomial
\begin{eqnarray*}
f(\lambda) & = & x^3 + \left({\frac {4096}{15}}\,{\lambda}^{3}-115\,\lambda
-146\,{\lambda}^{2}-{\frac {25}{24}}\right)\, x^2\\
&& \;  +  6\,{\lambda}^{2} \left( \lambda-1 \right)
 \left( 8192\,{\lambda}^{2}-7150\,\lambda-475 \right)\,x
+ 1080\,{\lambda}^{4} \left( \lambda-1 \right) ^{2}
\left( 2048\,\lambda-1805 \right).
\end{eqnarray*}

There are five obvious specializations with $\rho=20$
where we match our free parameter $\lambda$ (i.e.~the $I_0^*$ fiber) with another cusp.
This way, we obtain explicit equations for extremal elliptic
K3 surfaces over $\Q$
that were not derived in \cite{BM} or \cite{S-Rocky}.
These fibrations provide alternative realizations
of weight 3 newforms for the following discriminants
previously realized by semistable
extremal elliptic K3 surfaces (cf.~Section~\ref{s:Ex}):

\begin{center}
\begin{tabular}{|c|c|c|c|}
\toprule
$\lambda$ & degeneration in $\NS(X_\lambda)$
& disc $\NS(X_\lambda)$ & $T(X_\lambda)$\\
\midrule
\midrule
$-\frac{5}{1024}$ & $D_4 \rightsquigarrow D_5$ & $-840$ & [20, 0, 42]\\
\midrule
$\frac{35}{32}$ & $A_1 + D_4 \rightsquigarrow D_6$ & $-420$ & [6, 0, 70]\\
\midrule
$1$ & $A_2 + D_4 \rightsquigarrow D_7$ & $-280$ & [4, 0, 70]\\
\midrule
$0$ & $A_4 + D_4 \rightsquigarrow D_9$ & $-168$ & [4, 0, 42]\\
\midrule
$\infty$ & $A_6 + D_4 \rightsquigarrow D_{11}$ & $-120$ & [6, 0, 20]\\
\bottomrule
\end{tabular}
\end{center}

Every other specialization with $\rho=20$ requires
the existence of an additional section.
We ran Algorithm \ref{Alg} through the first 30 primes $p$
to find candidate parameters $\lambda_0\in\Q$.
In step (\ref{item}), we restricted to those $p$
which split in the field $K$ corresponding to a given newform $f$.
In almost every case,
this implied that there was exactly one candidate
parameter $\lambda\mod p$ (if there were any at all).
Some lifts to $\Q$ will be given in Table~\ref{Tab:fam}.
We shall now explain
how the structure of $X_\lambda$ as an elliptic surface
helps us verify that the K3 surface
at such a lift $\lambda_0$ is singular, i.e.~has $\rho=20$.

\subsection{Elliptic surfaces}
\label{ss:ell}

Once we know the group of sections of an elliptic surface,
we can easily compute the discriminant of its \NeS\ lattice using
the theory of \MoW\ lattices \cite{ShMW}.
On an elliptic K3 surface $X$,
the height pairing of two sections $P, Q$ is by \cite[Thm.~8.6]{ShMW}
\begin{eqnarray}
\label{eq:P,Q}
\langle  P, Q\rangle = 2 + (P.O) + (Q.O) - (P.Q)
- \sum_v \corr_v(P,Q).
\end{eqnarray}
The height pairing involves intersection numbers
in $\NS(X)$ and correction terms
according to the non-identity fiber components
met by the sections.
Since mostly we are concerned with
semistable elliptic fibers,
we give only the correction term at an $I_n$ fiber
(or root lattice $A_{n-1}$).
Using the cyclical numbering of components
such that $O$ meets $\Theta_0$, we have
\[
\corr_v(P, Q) = \frac{i\,(n-j)}n,\;\;\;
\text{if $P$ meets $\Theta_i$ and $Q$
meets $\Theta_j$ of a $I_n$ fiber}.
\]
If necessary, we interchange $P$ and $Q$ or renumber the components so that $i<j$.
In case $P=Q$, the height pairing
 specializes to the height $\hat h(P)$ of $P$:
 \begin{eqnarray}\label{eq:height}
\hat h(P) = 4 + 2\,(P.O) - \sum_v \corr_v(P).
\end{eqnarray}
The height pairing endows the \MoW\ group
modulo torsion with the structure of a positive-definite lattice,
the so-called \emph{\MoW\ lattice} $\MWL(X)$.
Note that the \MoW\ lattice
need not be integral if there are reducible fibers, neither does its discriminant.

One way to derive the height pairing is via
the orthogonal projection in $\NS(X)\otimes\Q$
with respect to the subspace generated by
zero section and fiber components.
Hence the discriminant of an elliptic surface $X$ satisfies
\begin{eqnarray}\label{eq:MWL}
(-1)^{\text{rk}(\MW(X))} \left|\MW(X)_\text{tor}\right|^2\,\disc(\NS(X)) =
\disc(\MWL(X)) \cdot \prod_v \disc(F_v).
\end{eqnarray}
Here the product runs over all reducible fibers $F_v$
and involves the discriminants of the corresponding root lattices.
For all elliptic fibrations in this paper,
the \MoW\ groups are always torsion-free
due to the configuration of singular fibers.

In the case at hand, the specialization $X$ at
$\lambda_0$ has a conjectural section $P\neq O$,
generating the \MoW\ group.
Specializing from \eqref{eq:P,Q} and \eqref{eq:MWL}
to the case of \MoW\ rank one,
 the discriminant is given by the product
of the discriminants of the root lattices and the height of $P$:
\begin{eqnarray}\label{eq:disc}
\disc\NS(X) = -2\cdot 3 \cdot 5\cdot 7 \cdot 4 \cdot \hat h(P)
= -840\, \hat h(P).
\end{eqnarray}
As explained above, the height of $P$ is determined by the
intersection number $(P.O)$ of $P$ and $O$
in $\NS(X)$ and the correction terms corr$_v(P)$
according to the fiber components which $P$ meets.
%
%
For the convenience of the reader, we record the possible correction terms in this situation:
\[
\corr_v(P)=\begin{cases}
0 &
\text{if $P$ meets the identity component of the fiber at $v$},\\
1 &
\text{if $P$ meets a non-identity component at a $I_0^*$ fiber},\\
\frac{i(n-i)}n &
\text{if $P$ meets the component $\Theta_i$ of a $I_n$ fiber}.
\end{cases}
\]

Assume now that we have a candidate parameter $\lambda_0$
for a specialization $X$ to realize some newform $f$.
Let $K$ denote the imaginary quadratic field corresponding to $f$.
Then we have to arrange for a section $P$\/ such that $\NS(X)$
has discriminant $d<0$ so that $\Q(\sqrt{d})=K$.
Let us look at one example in detail.

\begin{Example}\label{Ex:-1540}
Let $f$ be a newform for the field $K=\Q(\sqrt{-385})$
of discriminant $d=-1540$.
Then Algorithm \ref{Alg} suggests the candidate parameter
$\lambda_0 = 7^2\,11 / 2^9$.
We want to find a section $P$ such that disc $\NS(X)=d$.

Consider the formula (\ref{eq:disc}) for the discriminant
of the elliptic surface with section $P$.
Our target discriminant $d=-1540$ requires that
we eliminate a factor $6$ from (\ref{eq:disc})
while preserving $5$ and $7$.
Hence $P$\/ must meet the $I_5$ and $I_7$ fibers
on their identity components,
but meet the $I_2$ and $I_3$ fibers on non-identity components.
Now we assume that $(P.O)=0$ and that $P$ intersects
a non-identity component of the $I_0^*$ fiber.
Then  (\ref{eq:height}) reads
\[
\hat h(P) = 4 - \frac 12 - \frac 23 -1 = \frac{11}6.
\]
Hence (\ref{eq:disc}) gives
\[
\disc\NS(X) = -840 \cdot \frac{11}6 = -1540
\]
as required. The assumption $(P.O)=0$ implies
that $P$\/ has coordinates $(u, v)$ for polynomials $u, v$
of degree $4$ resp.~$6$ in the coordinate $t$ of the base curve.
Note that inserting $P$ into the Weierstrass equation of $X$
allows us to express the coefficients of $v$ in terms of those of $u$.
For the fibers to be met as prescribed, we must have
\[
(t-1)\left(t-\frac{35}{32}\right)(t-\lambda_0)^2 \mid v,\;\;\;
(t-1)(t-\lambda_0) \mid u.
\]
Moreover, the fibers at $\lambda_0$ and $\frac{35}{32}$
give two more linear relations in the coefficients of $u$.
Since $u$ has degree four, there is only one degree of freedom left.
Thus the problem is easily solved using a computer algebra system.
The solution $P$ and the resulting transcendental lattice
can be found in Table 1.
\end{Example}

\subsection{$p$-adic Newton iteration}
\label{ss:newton}

Often we are not lucky enough
to arrive at a system of equations that we can solve directly
as in Example~\ref{Ex:-1540}.
If we cannot find a direct solution, we apply the following algorithm,
essentially a $p$-adic Newton iteration in several variables.

\begin{Algorithm}\label{Newton}
Given a system $f_1=\hdots=f_n=0$
of algebraically independent polynomial equations
over $\Q$ in $n$ variables $z_1,\hdots,z_n$.
The following procedure tests for a solution over $\Q$:
\begin{enumerate}[(1)]
\item \label{Newton-1}
With an exhaustive search,
  find a solution $(\bar z_1,\hdots,\bar z_n)$ modulo some prime $p$.
\item \label{Newton-2}
Using difference quotients,
  double the $p$-adic accuracy of the solution a few times.
\item \label{Newton-3}
Compute a lift in $\Q$ with the Euclidean algorithm.
\item \label{Newton-4}
Test whether the lift solves the system of equations over $\Q$.
  While it does not, return to step (\ref{Newton-2}) to double the precision
  once more and try (\ref{Newton-3}) and (\ref{Newton-4}) again.
\end{enumerate}
\end{Algorithm}

For step~(\ref{Newton-2}) to converge,
we need some regularity assumptions for the polynomials~$f_i$.
For instance, it will converge if the coefficients of the $f_i$
are \hbox{$p$-adic} integers and the Jacobian determinant
$|\partial(f_1,\ldots,f_n) / \partial(z_1,\ldots,z_n)|$
does not vanish at $(\bar z_1,\hdots,\bar z_n)$.

We next outline some further implementation issues
specific to our setting.

Algorithm \ref{Newton} does not require that the section itself
be defined over $\Q$, only its $x$-coordinate.
Even if the $y$-coordinate involves
a square root of a rational number as a factor,
we can always arrange to solve a system of equations over $\Q$.



A more delicate point about Algorithm \ref{Newton} is
the choice of the prime $p$.
Here we distinguish whether or not $p$ splits
in the fixed imaginary quadratic field $K$.

For a singular K3 surface $X$ over a number field $L$,
one can predict the geometric Picard number
of the reductions.
Namely, let $\p$ denote a prime of $L$ above $p$.
If $X$ has good reduction mod $\p$,
write $X_\p$ for the reduced K3 surface. Then
\begin{eqnarray}\label{eq:Tate}
\rho(X_\p) =
\begin{cases}
20, & \text{if $p$ splits in $K$},\\
22, & \text{if $p$ is inert or ramified in $K$}.
\end{cases}
\end{eqnarray}
This follows from the \SI\ structure (cf.~sect.~\ref{s:obs})
since the above cases decide exactly
whether the elliptic curves $E, E'$
(and the abelian surface $E\times E'$) are supersingular.

In more detail, we apply the Tate conjecture \cite{Tate-C}
to a (conjectural) singular K3 surface $X$ over $\Q$.
Here we consider elliptic K3 surfaces with section,
so the Tate conjecture holds true by \cite{ASD}.
Note that the assumption implies
that the primes ramifying in $K$ are always bad,
since they divide the level of the associated Hecke eigenform
by Theorem~\ref{Thm:mod}.

At the inert primes $p$ in $K$,
the Fourier coefficient $a_p$ of the
associated newform is zero.
The resulting eigenvalues of Frob$_p$ are $p$ and $-p$.
Hence the Tate conjecture predicts
that the reduction $X_p$ has
an additional algebraic cycle over $\F_p$
and one more over $\F_{p^2}$.
On an elliptic surface, these extra cycles would either change the
configuration of reducible fibers or appear as extra \MoW\ sections.
In the former case we might miss the reduced surface entirely.
In the latter case, reduction mod~$p$ would increase
the \MoW\ rank compared to the rank over~$\Q$.
Depending on the conditions we impose ---
on $(P.O)$ and the fiber components met ---
step (\ref{Newton-1}) of Algorithm \ref{Newton} might then return more than
one section.
However, only one of these sections would lift to $\Q$,
so we would have to make the right choice.

Therefore, we always run Algorithm \ref{Newton}
at a prime $p$ that splits in $K$.
Under the assumption that the K3 surface is singular,
the Picard number in (\ref{eq:Tate}) guarantees
that the \MoW\ rank is constant at $20$
upon reduction.

\begin{Remark}\label{Rem:no-guess}
We will apply Algorithm \ref{Newton}
after guessing the parameter $\lambda_0$
with the help of Algorithm \ref{Alg}.
Hence there is one further equation $f_0=0$ in $z_1,\hdots,z_n$.
After step (\ref{Newton-4}), we thus also have to verify that $f_0$ vanishes at the lift.

The same ideas can also be applied
without a candidate parameter $\lambda_0$ at hand.
In Example~\ref{Ex:-1540}, this would have sufficed as well:
in the end, we would have to solve two equations in two variables.
However, if we have to apply Algorithm \ref{Newton}
to find the section explicitly,
it is computationally very convenient for step (\ref{Newton-1})
to have one parameter fewer.
\end{Remark}

\subsection{Specializations over $\Q$ with $\rho=20$}

Table \ref{Tab:fam} collects the specializations
of $X_\lambda$ where we verified $\rho=20$.
The additional section $P$ can be recovered
from its $x$-coordinate $u(t)$ by taking $x=u(t)$
in the Weierstrass equation and choosing a root for $v(t)$.
We also list the height $\hat h(P)$ of $P$
and the discriminant of $\NS(X)$.
The discriminant will be justified below.
The computation of the transcendental lattices $T(X)$
will be explained in the next subsection.

\begin{table}[ht!]
$$
\begin{array}{c|c|c|c|c}
\toprule
{\scriptstyle\lambda } & {\scriptstyle \text{disc } \NS(X) } & {\scriptstyle P: u(t) } & {\scriptstyle \hat h(P) } & {\scriptstyle T(X) } \\
\midrule
\midrule
{\scriptstyle -\frac 12} & {\scriptstyle -225\cdot 4} & 
{\scriptstyle -\frac{3^5}{2^3\cdot 5}\,(24\,t+35)\,(t-1)\,(2\,t+1)} &
{\scriptstyle \frac{15}{14}} & {\scriptstyle [30,0,30]}\\
\midrule
{\scriptstyle\frac 58 } & {\scriptstyle -4\cdot 43 } & {\scriptstyle
-{\frac {2\cdot 3^6}{5}}\,{t}^{2} \left( t-1 \right)
} & {\scriptstyle \frac{43}{210} } & {\scriptstyle [4,2,44]} \\
\midrule
{\scriptstyle\frac{5^3}{2^7} } & {\scriptstyle -4\cdot 67 } & {\scriptstyle
-{\frac {3^6}{2^6\cdot 5}}\, \left( t-1 \right) t \left( 128\,t-105 \right)
} & {\scriptstyle \frac{67}{210} } & {\scriptstyle [4,2,84]} \\
\midrule
{\scriptstyle\frac{5}{2^5} } & {\scriptstyle -88 } & {\scriptstyle
-{\frac {3^6}{2^4\cdot 5}}\, \left( t-1 \right) t \left( 32\,t-5 \right)
} & {\scriptstyle \frac{11}{105} } & {\scriptstyle [2,0,44]} \\
\midrule
{\scriptstyle\frac{5}{2^3\cdot 19^2} } & {\scriptstyle -4\cdot 163 } & {\scriptstyle
-{\frac {2\cdot 3^6\cdot 29^2}{5\cdot 19^2}}\,{t}^{2} \left(t-1 \right)
} & {\scriptstyle \frac{163}{210} } & {\scriptstyle [4,2,164]} \\
\midrule
{\scriptstyle\frac 52 } & {\scriptstyle -228 } & {\scriptstyle -\frac{3^6}{2^3\cdot 5}\,t\,(2\,t-5)\,(8\,t-7)
} & {\scriptstyle \frac{19}{70} } & {\scriptstyle [2,0,114]} \\
\midrule
{\scriptstyle\frac{3^2\cdot 5}{2^5} } & {\scriptstyle -312 } & {\scriptstyle
-{\frac {3^6}{2^4\cdot 5}}\,{t}^{2} \left( 32\,t-45 \right) } & {\scriptstyle \frac{13}{35}
} & {\scriptstyle [6,0,52]} \\
\midrule
{\scriptstyle-\frac 12 } & {\scriptstyle -340 } & {\scriptstyle
-{\frac {3^5}{2^3\cdot 5}}\, \left( t-1 \right)
\left( 2\,t+1 \right)  \left( 24\,t+35 \right)
} & {\scriptstyle \frac{17}{42} } & {\scriptstyle [20,10,22]} \\
\midrule
{\scriptstyle\frac{5\cdot 11^2}{2^9} } & {\scriptstyle -372 } & {\scriptstyle
-{\frac {3^6}{2^8\cdot 5^2}}\,{t}^{2} \left( 512\,t-605 \right) } & {\scriptstyle \frac{31}{70}
} & {\scriptstyle [6,0,62]} \\
\midrule
{\scriptstyle\frac{2\cdot 5}{3^2} } & {\scriptstyle -408 } & {\scriptstyle
\frac{3^3}{2^4\cdot 5}\, t\, (9\, t-10) (-96\, t+81)
} & {\scriptstyle \frac{17}{35} } & {\scriptstyle [6,0,68]} \\
\midrule
{\scriptstyle-\frac 5{2^3\cdot 3} } & {\scriptstyle -4\cdot 435 } & {\scriptstyle
\frac 3{2^4\cdot 5}\, (-7776\, t^3+4816\, t^2+1660\, t+175)
} & {\scriptstyle \frac{29}{14} } & {\scriptstyle [20,10,92]} \\
\midrule
{\scriptstyle\frac{3\cdot 5\cdot 7}{2^7} } & {\scriptstyle -4\cdot 483 } & {\scriptstyle
\frac{3^2}{2^8\cdot 5\cdot 7^2}\,(-5242880\, t^2+9202816\, t-3988061)
} & {\scriptstyle \frac{23}{10} } & {\scriptstyle [4,2,484]} \\
\midrule
{\scriptstyle-\frac 5{7^2} } & {\scriptstyle -520 } & {\scriptstyle
-{\frac {3^6}{2^4\cdot 5\cdot 7^2}}\, \left( t-1 \right)
\left( 32\,t+5 \right)  \left( 49\,t+5 \right)
} & {\scriptstyle \frac{13}{21} } & {\scriptstyle [20,0,26]} \\
\midrule
{\scriptstyle\frac{5\cdot 7^2}{2^5\cdot 19} } & {\scriptstyle
-532 } & {\scriptstyle -{\frac {2\cdot 3^5}{5^2\cdot 19}}\,
\left( t-1 \right)\, {t}^{2} \left( 608\,t-245 \right)
} & {\scriptstyle \frac {19}{30} } & {\scriptstyle [4,2,134]} \\
\midrule
{\scriptstyle\frac{37}{2^3\cdot 5} } & {\scriptstyle
-4\cdot 555 } & {\scriptstyle \frac{3^6}{2^4\cdot 5^4}\, (-31648\, t^3+86320\, t^2-78300\, t+23625)
} & {\scriptstyle \frac{37}{14} } & {\scriptstyle [4,2,556]} \\
\midrule
{\scriptstyle-\frac{7}{2^7} } & {\scriptstyle
-4\cdot 595 } & {\scriptstyle \frac{3}{2^6\cdot 5^3}\, (t-1)\, (-1048576\, t^3-716160\, t^2-45225\, t-875)
} & {\scriptstyle \frac{17}6 } & {\scriptstyle [20,10,124]} \\
\midrule
{\scriptstyle-\frac 1{2^5} } & {\scriptstyle -660 } & {\scriptstyle
\frac{3^5}{2^{14}\cdot 5}\, (32\, t+ 1) (-3072\, t^2+2592\, t+105)
} & {\scriptstyle \frac{11}{14} } & {\scriptstyle [20,10,38]} \\
\midrule
{\scriptstyle\frac{5\cdot 11^2}{2\cdot 7^2} } & {\scriptstyle -708 } & {\scriptstyle
\frac{3^2}{2^3\cdot 5^2 \cdot 7^3}\, t\, (98\, t-605) (-53176\, t+46585)
} & {\scriptstyle \frac{59}{70} } & {\scriptstyle [6,0,118]} \\
\midrule
{\scriptstyle-\frac{19}{2^5\cdot 7^2} } & {\scriptstyle -760 } & {\scriptstyle
\frac{3^5}{2^4\cdot 5^3 \cdot 7^2}\, (t-1)\, (1568\, t+19) (-75\, t-1)
} & {\scriptstyle \frac{19}{21} } & {\scriptstyle [20,0,38]} \\
\midrule
{\scriptstyle-\frac{1}{2^3\cdot 3\cdot 5}} & {\scriptstyle-4\cdot 795} & 
{\scriptstyle -\frac{3}{2^4\cdot 5^4}\,U(t)/(15360\,t+7)^2}
& {\scriptstyle \frac{53}{14}} & {\scriptstyle [20, 10, 164]}\\
\midrule
{\scriptstyle\frac{5^3\cdot 7}{2^9} } & {\scriptstyle -1092 } & {\scriptstyle
\frac{3^6}{2^8\cdot 5^4 \cdot 7^4} \,t^2\, (512\, t-875) (-4096\, t+2765)
} & {\scriptstyle \frac{13}{10} } & {\scriptstyle [6,0,182]} \\
\midrule
{\scriptstyle\frac{3^2\cdot 5\cdot 11}{2^{10}} } & {\scriptstyle -1320 } & {\scriptstyle
\frac{3^6}{2^{30}\cdot 5}\, (1024\, t-495)\, (-2097152\, t^2+5225472\, t-1607445)
} & {\scriptstyle \frac{11}7 } & {\scriptstyle [4,0,330]} \\
\midrule
{\scriptstyle\frac{11^2}{2^5\cdot 3} } & {\scriptstyle -1380 } & {\scriptstyle
\frac{3}{2^{16}\cdot 5}\, (96\, t-121)\, (-4427776\, t^2+9923936\, t-5636785)
} & {\scriptstyle \frac{23}{14} } & {\scriptstyle [6,0,230]} \\
\midrule
{\scriptstyle\frac{3\cdot 5^3}{2^5} } & {\scriptstyle -1428 } & {\scriptstyle
\frac{3^3}{2^4\cdot 5^3\cdot 7^2} \,t\, (32\, t-375) (32\, t^2-99880 \, t+87500)
} & {\scriptstyle \frac{17}{10} } & {\scriptstyle [6,0,238]} \\
\midrule
{\scriptstyle\frac{7^2\cdot 11}{2^9} } & {\scriptstyle -1540 } & {\scriptstyle
-\frac{3}{2^8\cdot 5} \, (512\, t-539)\, (t-1)\, (512000\, t^2-1097257\, t+588245)
} & {\scriptstyle \frac{11}{6} } & {\scriptstyle [22,0,70]} \\
\midrule
{\scriptstyle-\frac{5\cdot 7}{2^{10}\cdot 3^2} } & {\scriptstyle -1848 } & {\scriptstyle
-\frac{3^4}{2^7\cdot 5\cdot 7^3} \, t\, (9216\, t+35)
(75582720\, t^2+631582\, t+1323)
} & {\scriptstyle \frac{11}5 } & {\scriptstyle [42,0,44]} \\
\bottomrule
\end{array}
$$
\caption{Singular specializations over $\Q$ in the family $X_\lambda$}
\label{Tab:fam}
\label{table}
\end{table}

For the largest discriminant $-4\cdot 795$, the section $P$ is not integral. Here $U(t)$, the numerator of $u(t)$, is given by
\begin{eqnarray*}
U(t) & = & 229323571200000\,t^5-191371714560000\,t^4-9553942361376\,t^3\\
&& \;\; -151103350160\,t^2-953437100\,t-2100875.
\end{eqnarray*}

We note that the table is not complete.
For each of the four discriminants
 $$
  d =-4 \cdot 1435,\;\; -4 \cdot 1155,\;\; -4 \cdot 1995,\;\; -5460, 
 $$
 we find using Algorithm~10 the lift
 $$
 \lambda_0 = 
 -\frac{5 \cdot 7 \cdot 11^2}{2^3\cdot 41},\;\;
 -\frac{7}{2^3},\;\;
 \frac{7\cdot 19}{2^7\cdot 5},\;\;
 \frac{3 \cdot 5 \cdot 7^2}{2 \cdot 19^2}
 $$
 respectively for which we expect $X_{\lambda_0}$ to be singular
 with \NeS\ discriminant~$d$\/; but the additional section is
 too complicated for us to compute easily, even using Algorithm~12.
 In the next section we realize each of the remaining nine discriminants,
 including these four, in a singular K3 surface not in this family $X_\lambda$.

\begin{Lemma}\label{Lem:NS}
For all K3 surfaces in the table,
$\NS(X)$ is generated by fiber components
and the sections $O$ and $P$.
In particular, the discriminant is as claimed.
\end{Lemma}

\emph{Proof:}
The fiber components together with the sections
$O$ and $P$ generate a lattice $N$ of rank 20.
We must show that $N=\NS(X)$.
The discriminant $d$ of $N$ is given
by (\ref{eq:disc}) and (\ref{eq:height}).
If the index of $N$ in $\NS(X)$ were greater than~$1$,
then $\NS(X)$ would have discriminant $d$
divided by the square of the index.

For any elliptic K3 surface $X_\lambda$ in our family,
the \MoW\ group is torsion-free.
Using the formula~(\ref{eq:height}), we see that
$2\hat h(P)$ is a \hbox{$2$-adic} integer for any section~$P$,
because $(P.O)$ is an integer and each correction term $\corr_v(P)$
is in $\frac12 \Z_2$.
By (\ref{eq:disc}), we derive the general relation
\begin{eqnarray}\label{eq:4}
4 \mid \disc(\NS(X_\lambda)).
\end{eqnarray}
In each case in the table,
the quotient $d/\!\disc(\NS(X))$ can therefore only be a square
if it equals $1$. The claim follows. \qed

\subsection{Transcendental lattices}

The transcendental lattices for the singular K3 surfaces above
were computed using lattice theory as developed by Nikulin \cite{N}.
Here we sketch the argument.

Given an even integral non-degenerate lattice $L$,
we denote its dual by $L^\vee$.
In \cite{N}, Nikulin introduced a quadratic form
on the quotient $L^\vee/L$ which he called the discriminant form:
\begin{eqnarray*}
q_L\0:\;\; L^\vee/L & \to & \Q \mod 2\Z\\
x & \mapsto & x^2 + 2\Z
\end{eqnarray*}

For each singular K3 surface $X$ in the table,
we know the \NeS\ lattice $\NS(X)$
by Lemma \ref{Lem:NS}.
Hence we can compute its discriminant form.

\begin{Theorem}[Nikulin {\cite[Prop.~1.6.1]{N}}]
Let $N$ be an even integral unimodular lattice.
Let $L$ be a primitive non-degenerate sublattice and $M=L^\bot$.
Then
\[
q_L\0 = -q_M\0.
\]
\end{Theorem}

Since $\NS(X)$ always embeds primitively into $H^2(X,\Z)$,
the theorem provides the discriminant form
of the transcendental lattice $T(X)$.

\begin{Theorem}[Nikulin {\cite[Cor.~1.9.4]{N}}]
The genus of an even integral lattice is determined
by its signature and discriminant form.
\end{Theorem}

We now use the fact that each singular K3 surface $X$
in the table is defined over $\Q$.
By Theorem~\ref{Thm:genus}, the genus of  the transcendental
lattice $T(X)$ consists of a single class.
Hence  the discriminant form of the \NeS\ lattice $\NS(X)$
determines $T(X)$ uniquely via the above two theorems.
Thus the computation is completed by verifying
that for the given transcendental lattices $q_{T(X)}\0=-q_{\NS(X)}\0$.

\section{The remaining discriminants}
\label{s:rem}

So far we have matched all but nine newforms of weight 3
from Section~\ref{s:CM}
with singular K3 surfaces over $\Q$.
For some of the remaining discriminants, we found candidate surfaces
in other one-dimensional families of K3 surfaces with
$\rho\geq 19$, see Examples~\ref{Ex:1}--\ref{Ex:4}.
For the other forms, we used slightly different techniques
to derive elliptic K3 surfaces designed for those particular forms.
These will be sketched in Examples~\ref{Ex:5}--\ref{Ex:8}.
The transcendental lattices are computed
using the discriminant form as before.

\subsection{Specializations in one-dimensional families}

\begin{Example}[Discriminants $-1155, -1995$]
\label{Ex:1}
In the family $X_{\lambda, \mu}$ we choose $\lambda$
to merge fibers of type $I_0^*$ and $I_1$.
 A general member $X$ of the resulting family has
\[
\NS(X) = U + A_2 + A_4 + A_6 + D_5,
\]
so $\rho(X)\geq 19$. For a specialization to be defined over $\Q$,
we furthermore need a rational cusp at an $I_1$ fiber,
i.e.~the cubic factor of $\Delta$ encoding the $I_1$ fibers
in terms of $\mu$ must have a rational zero.
This can be achieved by the following substitution:
\[
\mu=9\,{\frac { \left( \nu+1 \right) ^{3}}{5\,{\nu}^{3}+15\,{\nu}^{2}-5\,\nu+1}}.
\]
Then the rational cusp gives
\[
\lambda = -{\frac { \left( \nu-3 \right)  \left( 5\,{\nu}^{3}+
15\,{\nu}^{2}-5\,\nu+1 \right) }{ \left( 7\,{\nu}^{2}+1 \right) ^{2}}}
\]
Algorithm \ref{Alg} suggests several
singular specializations over $\Q$
which can be verified explicitly using Algorith \ref{Newton}
as sketched in \ref{ss:newton}.
Here we give two of the specializations.
The corresponding newforms seem to occur in the family
$X_\lambda$ from the previous section as well,
but there the conjectural sections have twice the same height
because of the relation (\ref{eq:4}),
so they would be much harder to find.


{\bf Discriminant $\mathbf{-1155}$:}
Let $\nu=-3/5$.
Then there is a section $P$ of height $\hat h(P)=11/4$.
Its \hbox{$x$-coordinate} is given by
\[
u(t)=-\frac{2\cdot 27}{5^8\cdot 7\cdot 11^5\cdot 13} \,
(242\, t-585)\, (46060586\, t^3+422472710\, t^2+32588325\, t+8292375).
\]
This singular K3 surface has discriminant $d=-1155$
and transcendental lattice $T(X)=[6,3,194]$.

{\bf Discriminant $\mathbf{-1995}$:}
Let $\nu=9/35$.
Then there is a section $P$ of height $\hat h(P)=19/4$.
Its \hbox{$x$-coordinate} is given by
\[
u(t)  =  -\frac{27\cdot 11^3}{2^5\cdot 5^{10}\cdot 7^7\cdot 53} \,
\dfrac{(784\, t-795)\,U(t)}{(8757\, t-9010)^2},
\]
where
\begin{eqnarray*}
U(t) & = & 519278509294553530368\, t^5-2767640394056706623700\, t^4\\
&& {} +5908183745712577772625\, t^3 -6312492415348218806875\, t^2\\
&& {} +3374618170228790821875\, t-721947602876973103125.
\end{eqnarray*}

This singular K3 surface has discriminant $d=-1995$
and transcendental lattice $T(X)=[46,11,46]$.
\end{Example}

For the remaining discriminants,
we constructed other suitable families of elliptic K3 surfaces
and applied the same techniques as before.
Here we give only the specializations in extended Weierstrass form
\begin{eqnarray}\label{eq:Weier}
X:\;\;\; y^2 \, = \, x^3 +A\,x^2 + B\,x+C.
\end{eqnarray}

\begin{Example}[Discriminant $-627$]
\label{Ex:2}
We consider an elliptic K3 surface with singular fibers
$I_3, I_6, I_{11}$ at $\frac 17, 0, \infty$:
\begin{eqnarray*}
A & = & \frac {25}{24}}\,{t}^{4}+
{\frac {293}{6}}\,{t}^{3}-{\frac {23645}{16}}\,{t}^{2}
\mbox{}+{\frac {1705}{12}}\,t-{\frac {1331}{384},\\
B & = & - \left( 7\,t-1 \right) {t}^{2} \left( 200\,{t}^{3}
+9276\,{t}^{2}-92442\,t+4477 \right),\\
C & = & 96\, \left( 7\,t-1 \right) ^{2}{t}^{4}
\left( 100\,{t}^{2}+4588\,t-15059 \right).
\end{eqnarray*}
It has a section $P$ of height $\hat h(P)=19/6$ and $x$-coordinate
\[
u(t) = -\frac{3}{2\cdot 7^2\cdot 11^5}\, t \,
(52734375\,{t}^{3}+538828125\,{t}^{2}-2025538427\,t+1004475087).
\]
The transcendental lattice is $T(X)=[22,11,34]$ with discriminant $d=-627$.
\end{Example}

\begin{Example}[Discriminant $-715$]
\label{Ex:3}
We consider an elliptic K3 surface with singular fibers
$I_4, I_5, I_{11}$ at $1, 0, \infty$:
\begin{eqnarray*}
A & = & -{\frac {19487171}{3808800}}\,{t}^{4}
+{\frac {3674891}{13800}}\,{t}^{3}
\mbox{}-{\frac {247797}{80}}\,{t}^{2}
+{\frac {37743}{8}}\,t+{\frac {23805}{32}},\\
B & = & -\frac 25\,t \left( t-1 \right)
\left( 161051\,{t}^{3}-5251521\,{t}^{2}+16877745\,t+8212725 \right),\\
C & = & -152352\,{t}^{2} \left( t-1 \right) ^{2}
\left( 1331\,{t}^{2}-17526\,t-23805 \right).
\end{eqnarray*}
It has a section $P$ of height $\hat h(P)=13/4$ and $x$-coordinate
\[
u(t) = {\frac {15}{11776}}\, \left( t-1 \right)
\left( 44289025\,{t}^{3}-35970275\,{t}^{2}+11995075\,t-1058529 \right).
\]
The transcendental lattice is $T(X)=[22,11,38]$
with discriminant $d=-715$.
\end{Example}

\begin{Example}[Discriminant $-1435$]
\label{Ex:4}
We consider an elliptic K3 surface
with singular fibers $I_5, I_7, I_{8}$ at $1, 0, \infty$:
\begin{eqnarray*}
A & = & -{\frac {16807}{332928}}\,{t}^{4}-{\frac {490}{2601}}\,{t}^{3}
+{\frac {333135767}{13872}}\,{t}^{2}
\mbox{}-{\frac {275656745}{5202}}\,t+{\frac {603351125}{20808}},\\
B & = & -{t}^{2} \left( t-1 \right)
\left( t+50 \right)  \left( 2401\,{t}^{2}-114044\,t+114244 \right),\\
C & = & -83232\,{t}^{4} \left( t-1 \right) ^{2}
\left( 343\,{t}^{2}+436\,t-3380 \right).
\end{eqnarray*}
It has a section $P$ of height $\hat h(P)=41/8$ and $x$-coordinate
\[
u(t) = {\frac {51^2}{2\cdot 5^4\cdot 7}}\,
{\frac {U(t)}{ \left( 20003760\,t+208409617 \right) ^{2}}},
\]
where
\begin{eqnarray*}
U(t) & = & 94791757788196875\,{t}^{5}
-13440531435036024375\,{t}^{4}\\
&& {} + 311827388703362736750\,{t}^{3}
-2250368299914898266350\,{t}^{2}\\
&& {} + 5701998864279209056695\,t
-4700672234454567466251 .
\end{eqnarray*}
The transcendental lattice is $T(X)=[38,3,38]$
with discriminant $d=-1435$.
\end{Example}

The four remaining discriminants
require more work for two reasons.
For three of them ($-1012, -3003, -3315$),
the discriminant has two large prime factors.
Hence representing one of them
by a singular fiber of an elliptic fibration
would be too much a restriction on the other singular fibers
in a one-dimensional family.
Instead, we use two-dimensional families
to find singular K3 surfaces over $\Q$
with these discriminants.
This approach will be sketched
in \ref{ss:2}.

\subsection{Discriminant $d=-5460$}
\label{ss:5460}

The discriminant of maximal absolute value, $d=-5460$, illustrates the constraints
that we are facing:
The Hilbert class field of $\Q(\sqrt{d})$
has Galois group $(\Z/2\Z)^4$.
Hence Theorem~\ref{Thm:NS} implies
that a family of K3 surfaces
with a specialization of discriminant $d$ over $\Q$
has to be fairly complicated.
Indeed families of elliptic K3 surfaces with big Picard number
and more than three semistable reducible fibers
tend to be hard to parametrize.

In fact, one could try to do with one-dimensional families.
For instance, Algorithm \ref{Alg} suggests
that both one-dimensional families
in Section~\ref{s:fam} and Example~\ref{Ex:1}
admit a specialization over $\Q$ with discriminant $N^2 d$ for some $N\in\N$. 
However, the conjectural section would have too big height
for explicit computations with Algorithm \ref{Newton}.

In other words, we must find the right balance
between two competing constraints:
a family of K3 surfaces
that is not too complicated to parametrize,
but admits the required Galois action;
and a specialization with the given discriminant
that is not too complicated to allow for
an explicit verification by Algorithm \ref{Newton}.

These complexity problems can be circumvented as follows.
We let go the first step of
guessing the candidate parameter for a singular specialization
on a family of K3 surfaces over $\Q$ by Algorithm \ref{Alg}.
Instead we apply Algorithm \ref{Newton} directly
to a family of elliptic K3 surfaces over $\F_p$.
To start the algorithm, we do not have to parametrize the family explicitly
since
we need only determine
all family members over $\F_p$ --- i.e.~a finite set of solutions to  a given system of equations.
By point counting, we can already filter some surfaces with the
Algorithm \ref{Alg} from Section~\ref{s:tech}.
After finding an appropriate section on one of these surfaces,
we increase the $p$-adic accuracy and finally compute a lift
for both the surface and the section (if existing).
Then we verify as before
that the lifted surface has the prescribed
configuration of singular fibers
and that the section lifts to this surface with prescribed height.

\begin{Example}[Discriminant $-5460$]
\label{Ex:5}
Using the above approach for $p=37$,
we found an elliptic K3 surface $X$ with the following singular fibers:
$$
\begin{array}{|c||c|c|c|c|c|c|}
\toprule
\text{cusp} & \infty & 0 & 1 & \frac{7^3\cdot 11}{23}
& \alpha, \alpha^\sigma & -\frac{1}{125}\\
\midrule
\midrule
\text{fiber} & I_7 & I_5 & I_4 & I_3 & I_2 & I_1\\
\bottomrule
\end{array}
$$
Here $\alpha, \alpha^\sigma$ are the roots of the
polynomial $7625\,t^2-1367158\,t-57967$.
The surface is given in extended Weierstrass form
(\ref{eq:Weier}) with coefficients
\begin{eqnarray*}
A & = & 3125\,{t}^{4}-784700\,{t}^{3}
-40778898\,{t}^{2}-18971036\,t-218491,\\
B & = & -2^{20}\cdot 3^5\,{t}^{2} \left( t-1 \right)
\left( 625\,{t}^{3}-151380\,{t}^{2}-1599171\,t+62426 \right),\\
C & = & 2^{38}\cdot 3^{10}\,{t}^{4} \left( t-1 \right) ^{2}
\left( 125\,{t}^{2}-29164\,t-17836 \right).
\end{eqnarray*}
It has a section $P$ of height $\hat h(P)=13/4$ and $x$-coordinate
\begin{eqnarray*}
u(t) & = & - \frac{2^{19}}{3^3\cdot 7^4}  \, \left( t-1 \right)
\left( 45273407\,{t}^{3}-3678666\,{t}^{2}+168432\,t-5324 \right).
\end{eqnarray*}
Details of the computation will be published in  \cite{E-5460}.
The transcendental lattice is $T(X)=[42,0,130]$
with discriminant $-5460$.
\end{Example}

\subsection{Specializations in two-dimensional families}
\label{ss:2}

For the remaining three discriminants,
we work with two-dimensional families
of K3 surfaces with $\rho\geq 18$
because of the two large prime factors of the discriminant.
We determine a specialization
with two independent sections $P, Q$.
The discriminant is then computed in terms of the
height pairing on the \MoW\ lattice \cite{ShMW}
as sketched in \ref{ss:ell}.

\begin{Example}[Discriminant $-1012$]\label{Ex:6}
We consider an elliptic K3 surface $X$ with singular fibers
$I_3, I_5, I_{11}$ at $t=3861/28124$, $0$, and~$\infty$:
\begin{eqnarray*}
A & = & 66\, \left( 351384\,{t}^{4}-372196\,{t}^{3}
+113098\,{t}^{2}-13539\,t+594 \right),\\
B & = & {\frac {2^7\cdot 11}{3}}\,t \left( 28124\,t-3861 \right)
\left( 47916\,{t}^{3}-44484\,{t}^{2}+9479\,t-594 \right) ,\\
C & = & 2^{11}\cdot 11\,{t}^{2} \left( 28124\,t-3861 \right) ^{2}
\left( 242\,{t}^{2}-193\,t+22 \right).
\end{eqnarray*}
The elliptic surface $X$ 
lives inside a two-dimensional family
which can be parametrized explicitly
(Example \ref{Ex:3} comes from a subfamily thereof).
The specialization $X$ has sections $P, P^\sigma$
that are conjugate over $\Q(\sqrt{23})$.
We found them by an extension of Algorithm \ref{Newton}
using lattice reduction instead of the Euclidean algorithm
to recognize algebraic numbers of degree greater than~$1$.
The sections have height $\hat h(P) = \hat h(P^\sigma) = 38/15$
and pairing $\langle P, P^\sigma\rangle=8/15$.
Let $$
U(t) =  ( 16427508+3425424\, \sqrt{23} ) {t}^{2}
+ (2159201+450164\, \sqrt{23}) t
\mbox{}+293522+61224\, \sqrt{23}.$$
The $x$-coordinates of the sections $P$\/ and $P^\sigma$ are
$$
u(t) \, = \, 3\,{\frac { \left( 6- \sqrt{23} \right) ^{2}t
\left( 28124\,t-3861 \right)   U(t) }{ \left( 392
\mbox{}+95\, \sqrt{23} \right) ^{2}}}
$$
%
%
and the Galois conjugate of~$u$.
By (\ref{eq:MWL}), the elliptic surface $X$ has discriminant
\begin{eqnarray*}\label{eq:2-dim}
\disc(\NS(X)) = - 3\cdot 5\cdot 11\cdot
\left|\begin{matrix}\langle P,P\rangle & \langle P, P^\sigma\rangle\\ \langle P, P^\sigma\rangle
& \langle P^\sigma, P^\sigma\rangle\end{matrix}\right| = -1012.
\end{eqnarray*}
We compute the transcendental lattice $T(X)=[22,0,46]$.
\end{Example}

For the final two discriminants,
we did not parametrize the two-dimensional families explicitly,
but rather employed the approach over $\F_p$
as explained in \ref{ss:5460} for Example~\ref{Ex:5}.

\begin{Example}[Discriminant $-3003$]\label{Ex:7}
We consider an elliptic K3 surface
in extended Weierstrass  form (\ref{eq:Weier}):
\begin{eqnarray*}
A & = & 2197\,t^4 - 48516\,t^3 + 393636\,t^2 - 1208764\,t + 411600,\\
B & = & 216\,t \, (3\,t-25) \, (10\,t-67) \,
(13351\,t^3-200168\,t^2+834514\,t-411600),\\
C & = & 108^2\,t^2 \, (3\,t-25)^2 \, (10\,t-67)^2 \,
(81133\,t^2-641164\,t+411600).
\end{eqnarray*}
This surface has the following reducible singular fibers:
$$
\begin{array}{|c||c|c|c|c|c|}
\toprule
\text{cusp} & \frac{49}{4} & \frac{25}{3} & \frac{67}{10} & 0 &  \infty\\
\midrule
\midrule
\text{fiber} & I_2 & I_3 & I_3 & I_6 & I_{7}\\
\bottomrule
\end{array}
$$
There are orthogonal sections $P,Q$\/ of heights
$\hat h(P)=11/6$ and $\hat h(Q)=13/6$.
Their $x$-coordinates are
\begin{eqnarray*}
P: & u(t) = 27 \, t \, (3\,t/25-1) \, (31671\,(t/25)^2 - 675\,t - 6700),\\
Q: & u(t) = -\frac{4}{21} \,t \, (5984\,t^3 - 65032\,t^2 + 16442\,t - 949725).
\end{eqnarray*}
By (\ref{eq:MWL}), $X$ has discriminant $-3003$.
The transcendental lattice is $T(X)=[6,3,502]$.
\end{Example}

\begin{Example}[Discriminant $-3315$]\label{Ex:8}
We consider an elliptic K3 surface $X$
in extended Weierstrass  form (\ref{eq:Weier}):
\begin{eqnarray*}
A & = & 1105 \, (152628125\,t^4-131340300\,t^3
+35302566\,t^2-4215388\,t+574821), \\
B & = &    - 3200 \cdot 11^3 \cdot 13 \cdot 17^2 \, t \, (119\,t-26)
\, (741\,t-289)\\
&& \;\;\;\;\; \times\, (1795625\,t^3-1311895\,t^2+317027\,t-33813), \\
C & = &    80^3 \cdot 11^6 \cdot 13^2\cdot 17^2 \, t^2 \,
(119\,t-26)^2 \, (741\,t-289)^2  \, (27625\,t^2-16594\,t+2601).
\end{eqnarray*}
This surface has the following reducible singular fibers:
$$
\begin{array}{|c||c|c|c|c|}
\toprule
\text{cusp} & \frac{26}{119} &  \frac{289}{741} &  0 &  \infty\\
\midrule
\midrule
\text{fiber} &I_3 &  I_4 &  I_5 &  I_{8}\\
\bottomrule
\end{array}
$$
There are orthogonal sections $P,Q$\/ of heights
$\hat h(P)=17/8$ and $\hat h(Q)=13/4$.  Their $x$-coordinates are
\begin{eqnarray*}
P: & u(t) = 44^3 \, (88179\,t^2 - 58882\,t + 2415744/221),\\
Q: & u(t) = -\frac{5}{208} \,  (741\,t-289)
      \, (1960038171\,t^3+503840883\,t^2-379180711\,t+5454297).
\end{eqnarray*}
By (\ref{eq:MWL}), $X$\/ has discriminant $-3315$.
The transcendental lattice is $T(X)=[2,1,1658]$.
\end{Example}

\section{Proof of Theorem \ref{thm}}
\label{s:proof}

Subject to ERH, there are 65 imaginary quadratic fields of class group exponent two
by \cite{Wb}.
By \cite{S-CM}, they are in bijection with newforms of weight 3 with rational eigenvalues
up to twisting.
All newforms of class numbers 1 and 2 can be realized by Kummer surfaces over $\Q$
by Lemma \ref{Lem:Kummer}.
For each other field of class group exponent 2 (class number 4, 8, or 16),
we have exhibited a singular K3 surface over $\Q$ in Sections \ref{s:obs}, \ref{s:fam} and \ref{s:rem}.
Each K3 surface comes with an elliptic fibration with section over $\Q$,
so we can consider a Weierstrass form
\[
y^2 = x^3 + Ax + B,\;\;\; A,B\in{\Q}[t].
\]
The $\Q(\sqrt{d})$-isomorphic models as in \eqref{eq:Kummer}
thus realize all newforms corresponding to the given imaginary quadratic field.
This completes the proof of Theorem \ref{thm}.

\subsection*{Acknowledgements} We thank K.~Hulek for suggesting
this problem. We benefitted greatly from 
many discussions with him and B.~van Geemen.
We thank the referee for his detailed comments
and J.~Cervi\~no for the pointer to Eichler's paper.
This collaboration began at the Clay Mathematics Institute 2006 Summer School
on  ``Arithmetic Geometry'' at the Mathemati\-sches Institut 
of Georg-August-Universit\"at in G\"ottingen.
Our thanks go to the organizers for the invitation.

\end{document}